\documentclass{article}
\usepackage{amssymb}

\usepackage{amsmath}

\input{tcilatex}

\begin{document}

\title{On the deformation theory of structure constants for associative
algebras}
\author{B.G. Konopelchenko \\
\\
Dipartimento di Fisica, Universita del Salento \\
and INFN, Sezione di Lecce, 73100 Lecce, Italy.e-mail:konopel@le.infn.it}
\maketitle

\begin{abstract}
Algebraic scheme for constructing deformations of structure constants for
associative algebras generated by a deformation driving algebras (DDAs) is
discussed. An ideal of left divisors of zero plays a central role in this
construction. Deformations of associative three-dimensional algebras with
the DDA being a three-dimensional Lie algebra and their connection with
integrable systems are studied.
\end{abstract}

\bigskip


 \bigskip

\section{Introduction}

An idea to study deformations of structure constants for associative
algebras goes back to the classical works of Gerstenhaber [1,2]. As one of
the approaches to deformation theory he suggested '' to take the point of
view that the objects being deformed are not merely algebras, but
essentially algebra with a fixed basis'' and to treat '' the algebraic set
of all structure constants as parameter space for deformation theory'' [2].

Thus, following this approach, one chooses the basis $\mathbf{P}_{0},%
\mathbf{P}_{1},...,\mathbf{P}_{N}$ for a given algebra\textit{\ A}, takes
the structure constants $C_{jk}^{n}$ defined by the multiplication table

\begin{equation}
\mathbf{P}_{j}\mathbf{P}_{k}=\sum_{n=0}^{N}C_{jk}^{n}\mathbf{P}_{n},\quad
j,k=0,1,...,N
\end{equation}
and look for their deformations $C_{jk}^{n}(x),$ where $%
(x)=(x^{1},...,x^{M}) $ is the set of deformation parameters, such
that the associativity condition

\begin{equation}
\sum_{m=0}^{N}C_{jk}^{m}(x)C_{ml}^{n}(x)=%
\sum_{m=0}^{N}C_{kl}^{m}(x)C_{jm}^{n}(x)
\end{equation}
or similar equation is satisfied.

A remarkable example of deformations of this type with M=N+1 has
been discovered by Witten [3] and Dijkgraaf-Verlinde-Verlinde [4].
They
demonstrated that the function F which defines the correlation functions $%
\langle \Phi _{j}\Phi _{k}\Phi _{l}\rangle =\frac{\partial ^{3}F}{\partial
x^{j}\partial x^{k}\partial x^{l}}$ \ etc in the deformed two-dimensional
topological field theory obeys the associativity equation (2) with the
structure constants given by

\begin{equation}
C_{jk}^{l}=\sum_{m=0}^{N}\eta ^{lm}\frac{\partial ^{3}F}{\partial
x^{j}\partial x^{k}\partial x^{m}}
\end{equation}
where $\ $constants $\eta ^{lm}=(g^{-1})^{lm}$ and $g_{lm}=\frac{\partial
^{3}F}{\partial x^{0}\partial x^{l}\partial x^{m}}$ where the variable $%
x^{0} $ is associated with the unite element. Each solution of the WDVV
equation (2), (3) describes a deformation of the structure constants of the
N+1- dimensional associative algebra of primary fields $\Phi _{j}$.

Interpretation and formalization of the WDVV equation in terms of Frobenius
manifolds proposed by Dubrovin [5,6] provides us with a method to describe
class of deformations of the so-called Frobenius algebras. An extension of
this approach to general algebras and corresponding F-manifolds has been
given by Hertling and Manin [7]. Beautiful and rich theory of Frobenius and
F-manifolds has various applications from the singularity theory to quantum
cohomology (see e.g. [6,8,9] ).

\ An alternative approach to the deformation theory of the structure
constants for commutative associative algebras has been proposed recently in
[10-14]. Within this method the deformations of the structure constants are
governed by the so-called central system (CS) . Its concrete form depends on
the class of deformations under consideration and CS contains, as particular
reductions, many integrable systems like WDVV equation, oriented
associativity equation, integrable dispersionless, dispersive and discrete
equations (Kadomtsev-Petviashvili equation etc). The common feature of the
coisotropic, quantum, discrete deformations considered in [10-14] is that
for all of them elements $p_{j}$ of the basis and deformation parameters $%
x_{j}$ form a certain algebra ( Poisson, Heisenberg etc). A general class of
deformations considered in [13] is characterized by the condition that the
ideal $J=<f_{jk}>$ generated by the elements $f_{jk}=-p_{j}p_{k}+%
\sum_{l=0}^{N}C_{jk}^{l}(x)p_{l}$ \ representing the multiplication table
(1) is closed. It was shown that this class contains a subclass of so-called
integrable deformations for which the CS has a simple and nice geometrical
meaning.

In the present paper we will discuss a pure algebraic formulation of such
integrable deformations. We will consider the case when the algebra
generating deformations of the structure constants, i.e. the algebra formed
by the elements $p_{j}$ of the basis and deformation parameters $x_{k}$(
deformation driving algebra (DDA)), is a Lie algebra. The basic idea is to
require that all elements $f_{jk}=-p_{j}p_{k}+%
\sum_{l=0}^{N}C_{jk}^{l}(x)p_{l}$ are left divisors of zero and
that they generate the ideal $J=<f_{jk}>$ of left divisors of
zero. This requirement gives rise to the central system which
governs deformations generated by DDA.

Here we will study the deformations of the structure constants for the
three-dimensional algebra in the case when the DDA is given by one of the
three-dimensional Lie algebras. Such deformations are parametrized by a
single deformation variable x . Depending on the choice of DDA and
identification of $p_{1},p_{2}$ and x with the elements of DDA, the
corresponding CS takes the form of the system of ordinary differential
equations or the system of discrete equations (multi-dimensional mappings).
In the first case the CS contains the third order ODEs from the Chazy-Bureau
list as the particular examples. This approach provides us also with the Lax
form of the above equations and their first integrals.

The paper is organized as follows. General formulation of the
deformation theory for the structure constants is presented in
section 2. Quantum, discrete and coisotropic deformations are
discussed in section 3. Three-dimensional Lie algebras as DDAs are
analyzed in section 4. Deformations generated by general DDAs are
studied in section 5. Deformations driven by the nilpotent and
solvable DDAs are considered in sections 6 and 7, respectively.

\section{Deformations of the structure constants generated by DDA}

So, we consider a finite-dimensional noncommutative algebra
\textit{A }with ( or without ) unite element $\mathbf{P}_{0}$. We
will restrict overself to a class of algebras which possess a
basis composed by pairwise commuting elements
$\mathbf{P}_{0},\mathbf{P}_{1},...,\mathbf{P}_{N}$. The table of
multiplication (1) defines the structure constants $C_{jk}^{l}$ .
The commutativity of the basis implies that $C_{jk}^{l}$
$=C_{kj}^{l}$. In the
presence of the unite element one has $C_{j0}^{l}=\delta _{j}^{l}$ where $%
\delta _{j}^{l}$ is the Kroneker symbol.

\ Following the Gerstenhaber's suggestion [1,2] we will treat the structure
constants $C_{jk}^{l}$ in a given basis as the objects to deform and will
denote the deformation parameters by $x^{1},x^{2},...,x^{M}$. \ For the
undeformed structure constants the associativity conditions (2) are nothing
else than the compatibility conditions for the table of multiplication (1).
In the construction of deformations we should first to specify a ''deformed
'' version of the multiplication table and then to require that this
realization is selfconsistence and meaningful.

\ Thus, to define deformations we

1) associate a set of elements $p_{0},p_{1},...,p_{N},x^{1},x^{2},...,x^{M}$
with the elements of the basis $\mathbf{P}_{0},\mathbf{P}_{1},...,\mathbf{P}%
_{N}$ and deformation parameters $x^{1},x^{2},...,x^{M}$,

2) consider the Lie algebra \textit{B} of the dimension N+M+1 with the basis
elements $e_{1},...,e_{N+M+1}$ obeying the commutation relations

\begin{equation*}
\left[ e_{\alpha },e_{\beta }\right] =\sum_{\gamma
=1}^{N+M+1}C_{\alpha \beta \gamma }e_{\gamma },\quad \alpha ,\beta
=1,2,...,N+M+1,
\end{equation*}

3) identify the elements $p_{0},p_{1},...,p_{N},x^{1},x^{2},...,x^{M}$ with
the elements $e_{1},...,e_{N+M+1}$ thus defining the deformation driving
algebra (DDA). Different identifications define different DDAs. We will
assume that the element $p_{0}$ is always a central element of DDA. The
commutativity of the basis in the algebra \textit{A }implies the
commutativity between $p_{j}$ and in this paper we assume the same property
for all $x^{k}$. So, we will consider the DDAs defined by the commutation
relations of the type

\begin{equation}
\left[ p_{j},p_{k}\right] =0, \left[ x^{j},x^{k}\right] =0, \left[
p_{0},p_{k}\right] =0,\left[ p_{0},x^{k}\right] =0, \left[ p_{j},x^{k}%
\right] =\sum_{l}\alpha _{jl}^{k}x^{l}+\sum_{l}\beta _{j}^{kl}p_{l}
\end{equation}
where $\alpha _{jl}^{k}$ and $\beta _{j}^{kl}$ are some constants,

4) consider the elements

\begin{equation*}
f_{jk}=-p_{j}p_{k}+\sum_{l=0}^{N}C_{jk}^{l}(x)p_{l},\quad
j,k=0,1,...,N
\end{equation*}%
of the universal enveloping algebra U(\textit{B}) of the algebra DDA(\textit{%
B}). These $f_{jk}$ ''represent'' the table (1) in U(\textit{B}). Note that $%
f_{j0}=f_{0j}=0$.

5) require the all $f_{jk}$ are non-zero left divisors of zero and have a
common right zero divisor.

In this case $f_{jk}$ generate the left ideal $J=<f_{jk}>$ of left
divisors of zero.
We remind that non-zero elements a and b are
called left and right divisors of zero if $ab=0$ (see e.g. [15]).

\textbf{\ Definition.} The structure constants $C_{jk}^{l}(x)$ are said to
define deformations of the algebra \textit{A} generated by given DDA if all $%
f_{jk}$ are left zero divisors with common right zero divisor.

To justify this definition we first observe that the simplest possible
realization of the multiplication table (1) in U(\textit{B}) given by the
equations $f_{jk}=0,j,k=0,1,...,N$ is too restrictive in general. Indeed,
for instance, for the Heisenberg algrebra B [12] such equations imply that $%
[p_{l},C_{jk}^{m}(x)]=0$ and , hence, no deformations are allowed. So, one
should look for a weaker realization of the multiplication table. A
condition that all $f_{jk}$ are just non-zero divisors of zero is a natural
candidate. Then, the condition of compatibility of the corresponding
equations $f_{jk}\cdot \Psi _{jk}=0,j,k=1,...,N$ where $\Psi _{jk}$ are
right zero divisors requires that the l.h.s. of these equations and, hence, $%
\Psi _{jk}$ should have a common divisor (see e.g. [15] ). We restrict
ourself to the case when $\Psi _{jk}=\Psi \cdot \Phi _{jk},j,k=1,...,N$
where $\Phi _{jk}$ are invertible elements of U(B). In this case one has the
compatible set of equations

\begin{equation}
f_{jk}\cdot \Psi =0,\quad j,k=0,1,...,N
\end{equation}
that is all left zero divisors $f_{jk}$ have common right zero
divisor $\Psi $.

These conditions impose constraints on $C_{jk}^{m}(x)$. To clarify these
constraints we will use the basic property of the algebra \textit{A}, i.e.
its associativity. First we observe that due to the relations (4) one has
the identity

\begin{equation*}
\lbrack p_{l},C_{jk}^{m}(x)]=\sum_{t=0}^{N}\Delta _{jk,l}^{mt}(x)p_{t}
\end{equation*}
where \ $\Delta _{jk,l}^{mt}(x)$ are certain functions of
$x^{1},...,x^{M}$ only. Then, taking into account (4), one obtains
the identity

\begin{equation}
(p_{j}p_{k})p_{l}-p_{j}(p_{k}p_{l})=\sum_{s,t=0}^{N}K_{klj}^{st}\cdot
f_{st}+\sum_{t=0}^{N}\Omega _{klj}^{t}(x)\cdot p_{t},\quad j,k,l=0,1,...,N
\end{equation}%
where

\bigskip
\begin{equation*}
K_{klj}^{st}=\frac{1}{2}(\delta _{k}^{s}\delta _{l}^{t}+\delta
_{k}^{t}\delta _{l}^{s})p_{j}-\frac{1}{2}(\delta _{k}^{s}\delta
_{j}^{t}+\delta _{k}^{t}\delta _{j}^{s})p_{l}+\frac{1}{2}(\delta
_{j}^{s}C_{kl}^{t}+\delta _{j}^{t}C_{kl}^{s})-
\nonumber\\
\frac{1}{2}(\delta _{l}^{s}C_{kj}^{t}+\delta
_{l}^{t}C_{kj}^{s})+\Delta _{kl,j}^{st}-\Delta _{kj,l}^{st}
\end{equation*}
and

\begin{equation*}
\Omega
_{klj}^{t}(x)=\sum_{s}C_{jk}^{s}C_{ls}^{t}-\sum_{s}C_{lk}^{s}C_{js}^{t}+%
\sum_{s,n}(\Delta _{kj,l}^{sn}-\Delta _{kl,j}^{sn})C_{sn}^{t}.
\end{equation*}

The identity (6) implies that for an associative algebra

\begin{equation}
\sum_{s,t=0}^{N}K_{klj}^{st}\cdot f_{st}+\sum_{t=0}^{N}\Omega
_{klj}^{t}(x)\cdot p_{t}=0,\quad j,k,l=0,1,...,N.
\end{equation}
Due to the relations (5) equations (7) imply that

\begin{equation*}
\left( \sum_{t=0}^{N}\Omega _{klj}^{t}(x)\cdot p_{t}\right) \Psi =0.
\end{equation*}
These equations are satisfied if

\begin{equation}
\Omega
_{klj}^{t}(x)=\sum_{s}C_{jk}^{s}C_{ls}^{t}-\sum_{s}C_{lk}^{s}C_{js}^{t}+%
\sum_{s,n}(\Delta _{kj,l}^{sn}-\Delta _{kl,j}^{sn})C_{sn}^{t}=0,\quad
j,k,l,t=0,1,..,N.
\end{equation}

This system of equations plays a central role in our approach. If
$\Psi $ has no left zero divisors linear in $p_{j}$ and U(B) has
no zero elements linear in $p_{j}$ then the relation (8) is the
necessary condition for existence of a common right zero divisor
for $f_{jk}$.

At N$\geq 3$ it is also a sufficient condition. Indeed, if
$C_{jk}^{m}(x)$ are such that equations (8) are satisfied then

\begin{equation}
\sum_{s,t=0}^{N}K_{klj}^{st}\cdot f_{st}=0,\quad j,k,l=0,1,...,N.
\end{equation}

Generically, it is the system of $\frac{1}{2}N^{2}(N-1)$ linear equations
for $\frac{N(N+1)}{2}$ unknowns $f_{st}$ with noncommuting coefficients $%
K_{klj}^{st}$. At $N\geq 3$ for generic (non zeros, non zero divisors) $%
K_{klj}^{st}(x,p)$ the system (9) implies that

\begin{equation}
\alpha _{jk}f_{jk}=\beta _{lm}f_{lm},\quad j,k,l,m=1,...,N
\end{equation}
and

\begin{equation}
\gamma _{jk}f_{jk}=0,\quad j,k=1,...,N
\end{equation}%
where $\alpha _{jk},\beta _{lm},\gamma _{jk}$ are certain elements
of U(B) ( see e.g. [16,17] ). Thus, all $f_{jk}$ are right zero
divisors. They are also left zero divisors. Indeed, due to Ado's
theorem ( see e.g. [18] ) finite-dimensional Lie algebra B and,
hence, U(B) are isomorphic to matrix algebras. For the matrix
algebras zero divisors ( matrices with vanishing determinants) are
both right and left zero divisors [15]. Then, under the assumption
that all $\alpha _{jk}$ and $\beta _{lm}$ are not zero divisors,
the relations (10) imply that the right divisor of one of $f_{jk}$
is also the right zero divisor for the others.

At N=2 one has only two relations of the type (10) and a right
zero divisor of one of $f_{11},f_{12},f_{22}$ is the right zero
divisor of the others. We note that it isn't that easy to control
assumptions mentioned above. Nevertheless, the equations (5) and
(8) certainly are fundamental one for the whole approach.

We shall refer to the system (8) as the Central System (CS) governing
deformations of the structure constants of the algebra \textit{A }generated
by a given DDA. Its concrete form depends strongly on the form of the
brackets $\left[ p_{t},C_{jk}^{l}(x)\right] $ which are defined by the
relations (4) for the elements of the basis of DDA. For stationary solutions
($\Delta _{jk,l}^{t}=0)$ the CS (8) is reduced to the associativity
conditions (2).

\section{Quantum, discrete and coisotropic deformations.}

\bigskip

Coisotropic, quantum and discrete deformations of associative
algebras considered in [10-14] represent particular realizations
of the above general scheme associated with different DDAs.

For the \textbf{quantum deformations }of noncommutative algebra
one has $M=N$ and the deformation driving algebra is given by the
Heisenberg algebra [12]. The elements of the basis of the algebra
\textit{A} and deformations parameters are identified with the
elements of the Heisenberg algebra in such a way that

\begin{equation}
\left[ p_{j},p_{0}\right] =0,\left[ p_{0},x^{k}\right] =0,\left[ p_{j},p_{k}%
\right] =0, \left[ x^{j},x^{k}\right] =0, \left[
p_{j},x^{k}\right] =\hbar \delta _{j}^{k}p_{0},\quad j,k=1,...,N
\end{equation}%
where $\hbar $ is the (Planck's) constant. For the Heisenberg DDA

\begin{equation}
\Delta _{jk,l}^{mt}=\hbar \delta _{0}^{t}\frac{\partial C_{jk}^{m}(x)}{%
\partial x^{l}}
\end{equation}%
and consequently

\begin{equation}
\Omega _{klj}^{n}(x)=\hbar \frac{\partial C_{jk}^{n}}{\partial x^{l}}-\hbar
\frac{\partial C_{kl}^{n}}{\partial x^{j}}%
+\sum_{m=0}^{N}(C_{jk}^{m}C_{ml}^{n}-C_{kl}^{m}C_{jm}^{n})=0,\quad
j,k,l,n=0,1,...,N.
\end{equation}

Quantum CS (14) governs deformations of structure constants for
associative algebra driven by the Heisenberg DDA. It has a simple
geometrical meaning of vanishing Riemann curvature tensor for
torsionless Christoffel symbols $\Gamma _{jk}^{l}$ identified with
the structure constants ($C_{jk}^{l}=\hbar \Gamma _{jk}^{l}$)
[12].

\ In the representation of the Heisenberg algebra (12) by
operators acting in a linear space $\mathit{H}$ left divisors of
zero are realized by operators with nonempty kernel. The ideal $J$
is the left ideal generated by operators $f_{jk}$ which have
nontrivial common kernel or, equivalently, for which equations

\begin{equation}
f_{jk}\left| \Psi \right\rangle =0,\quad j,k=1,2,...,N
\end{equation}%
have nontrivial common solutions $\left| \Psi \right\rangle \subset \mathit{H%
}$. The compatibility condition for equations (15) is given by the CS (14).
The common kernel of the operators $f_{jk}$ form a \ subspace \textit{H}$%
_{\Gamma }$ in the linear space \textit{H}. So, in the approach under
consideration the multiplication table (1) is realized only on \textit{H}$%
_{\Gamma }$ , but not on the whole \textit{H}. Such type of realization of
the constraints is well-known in quantum theory as the Dirac's recipe for
quantization of the first-class constraints [19]. In quantum theory context
equations (15) serve to select the physical subspace in the whole Hilbert
space. Within the deformation theory one may refer to the subspace \textit{H}%
$_{\Gamma }$ as the ''structure constants'' subspace. In [12] the recipe
(15) was the starting point for construction of the quantum deformations.

Quantum CS (14) contains various classes of solutions which describe
different classes of deformations. An important subclass is given by
iso-associative deformations, i.e. by deformations for which the
associativity condition (2) is valid for all values of deformation
parameters. For such quantum deformations the structure constants should
obey the equations

\begin{equation}
\frac{\partial C_{jk}^{n}}{\partial x^{l}}-\frac{\partial C_{kl}^{n}}{%
\partial x^{j}}=0,\quad j,k,l,n=0,1,...,N.
\end{equation}
These equations imply that $C_{jk}^{n}=\frac{\partial ^{2}\Phi ^{n}}{%
\partial x^{j}\partial x^{k}}$ where $\Phi ^{n}$ are some functions while
the associativity condition (2) takes the form

\begin{equation}
\sum_{m=0}^{N}\frac{\partial ^{2}\Phi ^{m}}{\partial x^{j}\partial x^{k}}%
\frac{\partial ^{2}\Phi ^{n}}{\partial x^{m}\partial x^{l}}=\sum_{m=0}^{N}%
\frac{\partial ^{2}\Phi ^{m}}{\partial x^{l}\partial x^{k}}\frac{\partial
^{2}\Phi ^{n}}{\partial x^{m}\partial x^{j}}.
\end{equation}%
It is the oriented associativity equation introduced in [20,5]. Under the
gradient reduction $\Phi ^{n}=\sum_{l=0}^{N}\eta ^{nl}\frac{\partial F}{%
\partial x^{l}}$ equation (18) becomes the WDVV equation (2),(3).

\ Non iso-associative deformations for which the condition (16) is not valid
are of interest too. They are described by some well-known integrable
soliton equations [12]. In particular, there are the Boussinesq equation
among them for N=2 and the Kadomtsev-Petviashvili (KP) hierarchy for the
infinite-dimensional algebra of polynomials in the Faa' de Bruno basis [12].
In the latter case the deformed structure constants are given by

\begin{equation}
C_{jk}^{l}=\delta _{j+k}^{l}+H_{j-l}^{k}+H_{k-l}^{j},\quad j,k,l=0,1,2...
\end{equation}
with

\begin{equation}
H_{k}^{j}=\frac{1}{\hbar }\text{P}_{k}\left( -\hbar \widetilde{\partial }%
\right) \frac{\partial \log \tau }{\partial x^{j}},\quad j,k=1,2,3,...
\end{equation}
where $\tau $ is the famous tau-function for the KP hierarchy and $%
P_{k}\left( -\hbar \widetilde{\partial }\right) \doteqdot P_{k}\left( -\hbar
\frac{\partial }{\partial x^{1}},-\frac{1}{2}\hbar \frac{\partial }{\partial
x^{2}},-\frac{1}{3}\hbar \frac{\partial }{\partial x^{3}},...\right) $ where
$P_{k}\left( t_{1},t_{2},t_{3},...\right) $ are Schur polynomials defined by
the generating formula $\exp \left( \sum_{k=1}^{\infty }\lambda
^{k}t_{k}\right) =\sum_{k=0}^{\infty }\lambda ^{k}P_{k}\left( \mathbf{t}%
\right) .$

\textbf{Discrete deformations of }noncommutative associative
algebras are generated by the DDA with $M=N$ and commutation
relations

\begin{equation}
\left[ p_{j},p_{k}\right] =0,\quad \left[ x^{j},x^{k}\right] =0,\quad \left[
p_{j},x^{k}\right] =\delta _{j}^{k}p_{j},\quad j,k=1,...,N.
\end{equation}%
In this case

\begin{equation}
\Delta _{jk,l}^{mt}=\delta _{l}^{t}(T_{l}-1)C_{jk}^{m}(x),\quad
j,k,l,m,t=0,1,2,...,N
\end{equation}%
where for an arbitrary function $\varphi (x)$ the action of $T_{j}$ is
defined by $T_{j}\varphi (x^{0},...,x^{j},...,x^{N})=\varphi
(x^{0},...,x^{j}+1,....,x^{N}).$ The corresponding CS is of the form

\begin{equation}
C_{l}T_{l}C_{j}-C_{j}T_{j}C_{l}=0,\quad j,l=0,1,...,N
\end{equation}%
where the matrices $C_{j}$ are defined as $\left( C_{j}\right)
_{k}^{l}=C_{jk}^{l},j,k,l=0,1,...,N.$ The discrete CS (22) governs discrete
deformations of associative algebras.\  The CS (22)  contains, as particular
cases, the discrete versions of the oriented associativity equation, WDVV\
equation, Boussinesq equation, discrete KP hierarchy and Hirota-Miwa
bilinear equations for KP $\tau $-function.

\ For \textbf{coisotropic deformations of }commutative algebras [10,11] again%
\textbf{\ }$M=N$, but the DDA is the Poisson algebra with $p_{j}$ and $x^{k}$
identified with the Darboux coordinates, i.e.

\begin{equation}
\{p_{j},p_{k}\}=0,\quad \{x^{j},x^{k}\}=0,\quad \{p_{j},x^{k}\}=-\delta
_{j}^{k},\quad j,k=0,1,...,N.
\end{equation}
where $\{,\}$ is the standard Poisson bracket. The algebra
U(\textit{B}) is the commutative ring of functions and divisors of
zero are realized by functions with zeros. So, the functions
$f_{jk}$ should be functions with common set $\Gamma $ of zeros.
Thus, in the coisotropic case the multiplication table (1) is
realized by the set of equations [10]

\begin{equation}
f_{jk}=0,\quad j,k=0,1,2,...,N.
\end{equation}
Well-known compatibility conditon for these equations is

\begin{equation}
\left\{ f_{jk},f_{nl}\right\} \mid _{\Gamma }=0,\quad j,k,l,n=1,2...,N.
\end{equation}
The set $\Gamma $ is the coisotropic submanifold in $R^{2(N+1)}$.
The condition (25) gives rise to the following system of equations
for the structure constants

\begin{eqnarray}
\left[ C,C\right] _{jklr}^{m} &\doteqdot &\sum_{s=0}^{N}(C_{sj}^{m}\frac{%
\partial C_{lr}^{s}}{\partial x^{k}}+C_{sk}^{m}\frac{\partial C_{lr}^{s}}{%
\partial x^{j}}-C_{sr}^{m}\frac{\partial C_{jk}^{s}}{\partial x^{l}}%
\allowbreak -C_{sl}^{m}\frac{\partial C_{jk}^{s}}{\partial x^{r}}+  \notag \\
+C_{lr}^{s}\frac{\partial C_{jk}^{m}}{\partial x^{s}}-C_{jk}^{s}\frac{%
\partial C_{lr}^{m}}{\partial x^{s}}) =0
\end{eqnarray}
while the equations $\Omega _{klj}^{n}(x)=0$ have the form of associativity
conditions (2)

\begin{equation}
\Omega
_{klj}^{n}(x)=%
\sum_{m=0}^{N}(C_{jk}^{m}(x)C_{ml}^{n}(x)-C_{kl}^{m}(x)C_{jm}^{n}(x))=0.
\end{equation}

Equations (26) and (27) form the CS for coisotropic deformations [10]. In
this case $C_{jk}^{l}$ is transformed as the tensor of the type (1,2) under
the general tranformations of coordinates $x^{j}$ and the whole CS (26),
(27) is invariant under these tranformations [14]. The bracket $\left[ C,C%
\right] _{jklr}^{m}$ has appeared for the first time in the paper [21] where
the co-called differential concomitants were studied. It was shown in [18]
that this bracket is a tensor only if the tensor $C_{jk}^{l}$ obeys the
algebraic constraint (27). In the paper [7] the CS (26), (27) has appeared
implicitly as the system of equations which characterizes the structure
constants for F-manifolds. In [10] it has been derived as the CS governing
the coisotropic deformations of associative algebras.

The CS (26), (27) contains the oriented associativity equation, the WDVV
equation, dispersionless KP hierarchy and equations from the genus zero
universal Whitham hierarchy as the particular cases [10,11]. Yano manifolds
and Yano algebroids associated with the CS (26),(27) are studied in [14].

We would like to emphasize that for all deformations considered above the
stationary solutions of the CSs obey the global associativity condition (2).

\bigskip

\section{Three-dimensional Lie algebras as DDA.}

\bigskip

In the rest of the paper we will study deformations of associative algebras
generated by three-dimensional real Lie algebra $\mathit{L}$ . The complete
list of such algebras contains 9 algebras (see e.g. [18]). Denoting the
basis elements by $e_{1},e_{2},e_{3}$, \ one has the following nonequivalent
cases:

1) abelian algebra $\mathit{L}_{1}$,

2) general algebra $\mathit{L}_{2}$: $\left[ e_{1},e_{2}\right] =e_{1},\left[
e_{2},e_{3}\right] =0,\left[ e_{3},e_{1}\right] =0$,

3) nilpotent algebra $\mathit{L}_{3}:\left[ e_{1},e_{2}\right]
=0,[e_{2},e_{3}]=e_{1},\left[ e_{3},e_{1}\right] =0$,

4)-7) four nonequivalent solvable algebras : $\left[ e_{1},e_{2}\right] =0,%
\left[ e_{2},e_{3}\right] =\alpha e_{1}+\beta e_{2},\left[ e_{3},e_{1}\right]
=\gamma e_{1}+\delta e_{2}$ with $\alpha \delta -\beta \gamma \neq 0$,

8)-9) simple algebras $\mathit{L }_{8}=$ so(3) and $\mathit{L} _{9}=$so(2,1).

\ In virtue of the one to one correspondence between the elements of the
basis in DDA and the elements $p_{j}$, $x^{k}$ an algebra $\mathit{L} $
should has an abelian subalgebra and only one its element may play a role of
the deformation parameter $x$. For the original algebra \textit{A }and the
algebra \textit{B} one has two options:

1) \textit{A} is a two-dimensional algebra without unite element and \textit{%
B }=\textit{L}

2) \textit{A} is a three-dimensional algebra with the unite element and
\textit{B} = $\mathit{L}_{0}\oplus \mathit{L}$ where $\mathit{L }_{0}$ is
the algebra generated by the unite element $p_{0}$.

After the choice of \textit{B} one should establish a correspondence between
$p_{1},p_{2},x$ and $e_{1},e_{2},e_{3}$ defining DDA. For each algebra $%
\mathit{L}_{k}$ there are obviously, in general, six possible
identifications if one avoids linear superpositions. Some of them are
equivalent. The incomplete list of nonequivalent identifications is:

1) algebra $\mathit{L}_{1}:p_{1}=e_{1},p_{2}=e_{2},x=e_{3};$ DDA is the
commutative algebra with

\begin{equation}
\lbrack p_{1},p_{2}]=0,\quad [p_{1},x]=0,\quad [p_{2},x]=0.
\end{equation}

2) algebra $\mathit{L}_{2}$ :

case a) $p_{1}=-e_{2},p_{2}=e_{3},x=e_{1}$ ; the corresponding DDA is the
algebra $\mathit{L}_{2a}$ with the commutation relations

\begin{equation}
\lbrack p_{1},p_{2}]=0,\quad [p_{1},x]=x,\quad [p_{2},x]=0,
\end{equation}

case b) $p_{1}=e_{1},p_{2}=e_{3},x=e_{2}$ ; the corresponding DDA $\mathit{L}%
_{2b}$ is defined by

\begin{equation}
\lbrack p_{1},p_{2}]=0,\quad [p_{1},x]=p_{1},\quad [p_{2},x]=0.
\end{equation}

3) algebra $\mathit{L}_{3}$: $p_{1}=e_{1},p_{2}=e_{2},x=e_{3}$; DDA $\mathit{%
L }_{3}$ is

\begin{equation}
\lbrack p_{1},p_{2}]=0,\quad [p_{1},x]=0,\quad [p_{2},x]=p_{1},
\end{equation}

4) solvable algebra $\mathit{L}_{4}$ with $\alpha =0,\beta =1,\gamma
=-1,\delta =0:p_{1}=e_{1},p_{2}=e_{2},x=e_{3}$ ; DDA $\mathit{L}_{4} $ with

\begin{equation}
\lbrack p_{1},p_{2}]=0,\quad [p_{1},x]=p_{1},\quad [p_{2},x]=p_{2},
\end{equation}

5) solvable algebra $\mathit{L}_{5}$ at $\alpha =1,\beta =0,\gamma =0,\delta
=1:p_{1}=e_{1},p_{2}=e_{2},x=e_{3}$ ; DDA $\mathit{L} _{5}$ is

\begin{equation}
\lbrack p_{1},p_{2}]=0,\quad [p_{1},x]=p_{1},\quad [p_{2},x]=-p_{2}.
\end{equation}

For the second choice of the algebra $\mathit{B=\mathit{L} }_{0}\oplus
\mathit{L} $ mentioned above the table of multiplication \ (1) consists from
the trivial part $\mathbf{P}_{0}\mathbf{P}_{j}=\mathbf{P}_{j}\mathbf{P}_{0}=%
\mathbf{P}_{j},j=0,1,2$ \ and the nontrivial part

\begin{eqnarray}
\mathbf{P}_{1}^{2} =A\mathbf{P}_{0}+B\mathbf{P}_{1}+C\mathbf{P}_{2},  \notag
\\
\mathbf{P}_{1}\mathbf{P}_{2} =D\mathbf{P}_{0}+E\mathbf{P}_{1}+G\mathbf{P}%
_{2}, \\
\mathbf{P}_{2}^{2} =L\mathbf{P}_{0}+M\mathbf{P}_{1}+N\mathbf{P}_{2}.  \notag
\end{eqnarray}

For the first choice $\mathit{B=}\mathit{L} $ the multiplication table is
given by (34) with A=D=L=0.

It is convenient also to arrange the structure constants A,B,...,N into the
matrices $C_{1},C_{2}$ defined by $(C_{j})_{k}^{l}=C_{jk}^{l}$. One has

\begin{equation}
C_{1}=\left(
\begin{array}{ccc}
0 & A & D \\
1 & B & E \\
0 & C & G%
\end{array}
\right) ,\quad C_{2}=\left(
\begin{array}{ccc}
0 & D & L \\
0 & E & M \\
1 & G & N%
\end{array}%
\right) .
\end{equation}

In terms of these matrices the associativity conditions (2) are written as

\begin{equation}
C_{1}C_{2}=C_{2}C_{1}.
\end{equation}

\section{\textbf{Deformations generated by general DDAs}}

\bigskip

1. Commutative DDA (28) obviously does not generate any deformation. So, we
begin with the three-dimensional commutative algebra \textit{A }and \textbf{%
DDA }$\mathit{L} _{2a}$ defined by the commutation relations (29). These
relations imply that for an arbitrary function $\varphi (x)$

\begin{equation}
\lbrack p_{j},\varphi (x)]=\Delta _{j}\varphi (x),\quad j=1,2
\end{equation}
where $\ \Delta _{1}=x\frac{\partial }{\partial x},\Delta _{2}=0$.
Consequently, one has the following CS

\begin{equation}
\Omega _{klj}^{n}(x)=\Delta _{l}C_{jk}^{n}-\Delta
_{j}C_{kl}^{n}+\sum_{m=0}^{2}(C_{jk}^{m}C_{lm}^{n}-C_{kl}^{m}C_{jm}^{n})=0,%
\quad j,k,l,n=0,1,2.
\end{equation}
In terms of the matrices $C_{1}$ and $C_{2}$ defined above this CS has a
form of the Lax equation

\begin{equation}
x\frac{\partial C_{2}}{\partial x}=[C_{2},C_{1}].
\end{equation}

The CS (39) has all remarkable standard properties of the Lax equations (see
e.g. [20,21]): it has three independent first integrals

\begin{equation}
I_{1}=trC_{2},\quad I_{2}=\frac{1}{2}tr(C_{2})^{2},\quad I_{3}=\frac{1}{3}%
tr(C_{2})^{3}
\end{equation}
and it is equivalent to the compatibility condition of the linear problems

\begin{eqnarray}
C_{2}\Phi  &=&\lambda \Phi ,  \notag \\
x\frac{\partial \Phi }{\partial x} &=&-C_{1}\Phi
\end{eqnarray}%
where $\Phi $ is the column with three components and $\lambda $ is a
spectral parameter. Though the evolution in x described by the second linear
problem (41) is too simple, nevertheless the CS (38) or (39) have the
meaning of the iso-spectral deformations of the matrix $C_{2}$ that is
typical to the class of integrable systems (see e.g. [22,23]).

CS (39) is the system of six equations for the structure constants
D,E,G,L,M,N with free A,B,C:

\begin{eqnarray}
D^{\prime} =DB+LC-AE-DG,  \notag \\
L^{\prime} =DE+LG-AM-DN,  \notag \\
E^{\prime} =MC-EG-D,  \notag \\
M^{\prime} =E^{2}+MG-BM-EN-L, \\
G^{\prime} =GB+NC-CE-G^{2}+A,  \notag \\
N^{\prime} =GE-CM+D  \notag
\end{eqnarray}
where $D^{\prime} =x\frac{\partial D}{\partial x}$ etc. Here we will
consider only simple particular cases of the CS (42). First corresponds to
the constraint A=0, B=0, C=0, i.e. to the nilpotent $\mathbf{P}_{1}$. The
corresponding solution is

\begin{eqnarray}
D =\frac{\beta }{\ln x},\quad E=-\beta +\frac{\gamma }{\ln x},\quad G=\frac{1%
}{\ln x},L=\alpha \beta +2\beta ^{2}+\delta \ln x-\frac{\beta \gamma }{\ln x}%
,  \notag \\
M =\alpha \gamma +3\beta \gamma +\mu \ln x-\delta (\ln x)^{2}-\frac{\gamma
^{2}}{\ln x},\quad N=\alpha +\beta -\frac{\gamma }{\ln x}
\end{eqnarray}
where $\alpha ,\beta ,\gamma ,\delta ,\mu $ are arbitrary constants. The
three integrals for this solution are

\begin{equation}
I_{1}=\alpha ,I_{2}=\frac{1}{2}\alpha ^{2}+3\beta ^{2}+2\alpha \beta +\mu
,I_{3}=\frac{1}{3}((\alpha +\beta )^{3}-\beta ^{3})+(\alpha +\beta )(\mu
+\beta (\alpha +2\beta ))-\gamma \delta.
\end{equation}

The second example is given by the constraint B=0, C=1, G=0 for which the
quantum CS (14) is equivalent to the Boussinesq equation [12]. Under this
constraint the CS (42) is reduced to the single equation

\begin{equation}
E^{\prime \prime} -6E^{2}+4\alpha E+\beta =0
\end{equation}
and the other structure constants are given by

\begin{eqnarray}
A =2E-\alpha ,\quad B=0,\quad C=1,\quad D=\gamma -\frac{1}{2}E^{\prime}
,\quad G=0,  \notag \\
L =-E^{2}+\alpha E+\frac{1}{2}\beta ,\quad M=\gamma +\frac{1}{2}E^{\prime}
,N=\alpha -N
\end{eqnarray}
where $\alpha ,\beta ,\gamma $ are arbitrary constants. The corresponding
first integrals are

\begin{equation}
I_{1}=\alpha ,I_{2}=\frac{1}{2}(\beta +\alpha ^{2}),I_{3}=\frac{1}{3}\alpha
^{3}+\gamma ^{2}+\frac{1}{2}\alpha \beta -\frac{1}{4}(E^{\prime
})^{2}+E^{3}-\alpha E^{2}-\frac{1}{2}\beta E.
\end{equation}%
Integral $I_{3}$ reproduces the well-known first integral of equation (45).
Solutions of equation (45) are given by elliptic integrals (see e.g. [24]).
Any such solution together with the \ formulae (46) describes deformation of
the three-dimensional algebra \textit{A} driven by DDA $\mathit{\mathit{L}}%
_{2a}$.

Now we will consider deformations of the two-dimensional algebra \textit{A }%
without unite element according to the first option mentioned in the
previous section. In this case the CS has the form (39) with the $2\times 2 $
matrices

\begin{equation}
C_{1}=\left(
\begin{array}{cc}
B & E \\
C & G%
\end{array}%
\right) ,\quad C_{2}=\left(
\begin{array}{cc}
E & M \\
G & N%
\end{array}%
\right)
\end{equation}
or in components

\begin{eqnarray}
E^{\prime } =MC-EG,  \notag \\
M^{\prime } =E^{2}+MG-BM-EN,  \notag \\
G^{\prime } =GB+NC-CE-G^{2},  \notag \\
N^{\prime } =GE-CM.
\end{eqnarray}
In this case there are two independent integrals of motion

\begin{equation}
I_{1}=E+N,\quad I_{2}=\frac{1}{2}(E^{2}+N^{2}+2MG).
\end{equation}

The corresponding spectral problem is given by (41). Eigenvalues of the
matrix $C_{2}$, i.e. $\lambda _{1,2}=\frac{1}{2}(E+N\pm \sqrt{(E-N)^{2}+4GM%
\text{ }})$ are invariant under deformations and $\det C_{2}=\frac{1}{2}%
I_{1}^{2}-I_{2}$. We note also an obviously invariance of equations (42) and
(49) under the rescaling of x.

The system of equations (49) contains two arbitrary functions B and C. In
virtue of the possible rescaling $\mathbf{P}_{1}\rightarrow \mu _{1}\mathbf{P%
}_{1},\mathbf{P}_{2}\rightarrow \mu _{2}\mathbf{P}_{2}$ of the basis for the
algebra \textit{A }with two arbitrary functions $\mu _{1},\mu _{2}$, one has
four nonequivalent choices 1) B=0, C=0, 2) B=1, C=0, 3) B=0, C=1, 4) B=1,
C=1.

In the case B=0, C=0 ( nilpotent $\mathbf{P}_{1}$ ) the solution of the
system (49) is

\begin{equation}
B=0,C=0,E=\frac{\beta }{\ln x},G=\frac{1}{\ln x},M=\gamma \ln x-\frac{\beta
^{2}}{\ln x}+\alpha \beta ,N=-\frac{\beta }{\ln x}+\alpha
\end{equation}
where $\alpha $, $\beta ,\gamma $ are arbitrary constants. For this solution
the integrals are equal to $I_{1}=\alpha ,$ $I_{2}=\gamma +\frac{1}{2}\alpha
^{2}$ and $\lambda _{1,2}=\frac{1}{2}(\alpha +\sqrt{\alpha ^{2}+4\gamma })$.

At B=1, C=0 the system (49) has the following solution

\begin{eqnarray}
B=1,\quad C=0,\quad E=\frac{\gamma }{x+\beta },\quad G=\frac{x}{x+\beta },
\notag \\
M=\delta +(\alpha \gamma +\beta \delta -\frac{\gamma ^{2}}{\beta })\frac{1}{x%
}+\frac{\gamma ^{2}}{\beta (x+\beta )},\quad N=-\frac{\gamma }{x+\beta }%
+\alpha
\end{eqnarray}
where $\alpha ,\beta ,\gamma ,\delta $ are arbitrary constants. The
integrals are $I_{1}=\alpha ,$ $I_{2}=\delta +\frac{1}{2}\alpha ^{2}$. The
formulae (51), (52) provide us with explicit deformations of the structure
constants.

In the last two cases the CS (49) is equivalent to the simple third order
ordinary differential equations. At B=0, C=1 with additional constraint $%
I_{1}=0$ one gets

\begin{equation}
G^{\prime \prime \prime} +2G^{2}G^{\prime} +4(G^{\prime} )^{2}+2GG^{\prime
\prime} =0
\end{equation}
while at B=1,C=1 and $I_{1}=0$ the system (49) becomes

\begin{equation}
G^{\prime \prime \prime} +2G^{2}G^{\prime} +4(G^{\prime} )^{2}+2GG^{\prime
\prime} -G^{\prime} =0.
\end{equation}
The second integral for these ODEs is

\begin{equation}
I_{2}=-\frac{1}{2}G^{4}+\frac{1}{2}(G^{\prime })^{2}-2G^{2}G^{\prime
}-GG^{\prime \prime }+\frac{1}{2}BG^{2}.
\end{equation}%
Equation (53) with $G^{\prime }=\frac{\partial G}{\partial y}$ is the Chazy
V equation from the well-known Chazy-Bureau list of the third order ODEs
having Painleve property [25,26]. The integral (55) is known too (see e.g.
[27]).

The appearance of the Chazy V equation among the particular cases of the
system (49) indicates that for other choices of B and C the CS (49) may be
equivalent to the other notable third order ODEs. It is really the case.
Here we will consider only the reduction C=1 with $I_{1}=N+E=0$. In this
case the system (49) is reduced to the following equation

\begin{equation}
G^{\prime \prime \prime} +2G^{2}G^{\prime} +4(G^{\prime} )^{2}+2GG^{\prime
\prime} -2G^{\prime} \Phi -G\Phi^{ \prime} =0
\end{equation}
where $\Phi =B^{\prime} +\frac{1}{2}B^{2}$.The second integral is

\begin{equation}
I_{2}=-\frac{1}{2}G^{4}+\frac{1}{2}(G^{\prime} )^{2}-2G^{2}G^{\prime}
-GG^{\prime \prime} +\Phi G^{2}.
\end{equation}
and $\lambda _{1,2}=\pm \sqrt{\frac{I_{2}}{2}}$.

Choosing particular B or $\Phi $, one gets equations from the Chazy-Bureau
list. Indeed, at $\Phi =0$ one has the Chazy V equation (53). Choosing $\Phi
=G^{\prime} $, one gets the Chazy VII equation

\begin{equation}
G^{\prime \prime \prime} +2G^{2}G^{\prime} +2(G^{\prime} )^{2}+GG^{\prime
\prime} =0.
\end{equation}
At B=2G equation (56) becomes the Chazy VIII equation

\begin{equation}
G^{\prime \prime \prime} -6G^{2}G^{\prime} =0.
\end{equation}

Choosing the function $\Phi $ such that

\begin{equation}
\left( 6\Phi e^{\frac{1}{3}G}\right)^{ \prime} =2G^{2}G^{\prime} +(G^{\prime
})^{2}+4GG^{\prime \prime} ,
\end{equation}

one gets the Chazy III equation

\begin{equation}
G^{\prime \prime \prime }-2GG^{\prime \prime }+3(G^{\prime })^{2}=0.
\end{equation}%
In the above particular cases the integral $I_{2}$ (57) is reduced to those
given in [27].

All Chazy equations presented above have the Lax representation (39) with $%
E=-N=-\frac{1}{2}(G^{\prime} +G^{2}+GB),M=-\frac{1}{2}(G^{\prime \prime}
+3GG^{\prime} +G^{3}+G^{2}B+(GB)^{\prime }),C=1$ and the proper choice of B.

Solutions of all these Chazy equations provide us with the deformations of
the structure constants (48) for the two-dimensional algebra \textit{A}
generated by the DDA $\mathit{L} _{2a}$.

2.Now we pass to the \textbf{DDA} $\mathit{L} _{2b}.$ The commutation
relations (30) imply that

\begin{equation}
\lbrack p_{1},\varphi (x)]=(T-1)\varphi (x)\cdot p_{1},\quad
[p_{2}, \varphi (x)]=0
\end{equation}
where $\varphi (x)$ is an arbitrary function and $T\varphi (x)=\varphi
(x+1). $ Using (62), one finds the corresponding CS

\begin{eqnarray}
\sum_{m=0}^{2}((\Delta _{l}+1)C_{jk}^{m}(x)\cdot C_{lm}^{n}(x)=  \notag \\
=(\Delta _{j}+1)C_{kl}^{m}(x)\cdot C_{jm}^{n}(x)), \quad j,k,l,n=0,1,2
\end{eqnarray}
where $\Delta _{1}=T-1,\Delta _{2}=0.$ In terms of the matrices $C_{1}$ and $%
C_{2}$ this CS is

\begin{equation}
C_{1}TC_{2}=C_{2}C_{1}.
\end{equation}
For nondegenerate matrix $C_{1}$ one has

\begin{equation}
TC_{2}=C_{1}^{-1}C_{2}C_{1}.
\end{equation}

The CS (65) is the discrete version of the Lax equation (39) and has similar
properties. It has three independent first integrals

\begin{equation}
I_{1}=trC_{2},\quad I_{2}=\frac{1}{2}tr(C_{2})^{2},\quad I_{3}=\frac{1}{3}%
tr(C_{2})^{3}
\end{equation}
and represents itself the compatibility condition for the linear problems

\begin{eqnarray}
\Phi C_{2} =\lambda \Phi ,  \notag \\
T\Phi =\Phi C_{1}.
\end{eqnarray}
Note that $\det C_{2}$ is the first integral too.

The CS (64) is the discrete dynamical system in the space of the structure
constants. For the two-dimensional algebra \textit{A} with matrices (48) it
is

\begin{eqnarray}
BTE+ETG =EB+MC,  \notag \\
BTM+ETN =E^{2}+MG,  \notag \\
CTE+GTG =BG+CN, \\
CTM+GTN =EG+NG  \notag
\end{eqnarray}
where B and C are arbitrary functions. For nondegenerate matrix $C_{1},$
i.e. at $BG-CE\neq 0$ , one has the resolved form (65), i.e.

\begin{eqnarray}
TE =\frac{GM-EN}{BG-CE}C,\quad TG=B+\frac{BN-CM}{BG-CE}C,  \notag \\
TM =\frac{GM-EN}{BG-CE}G,\quad TN=E+\frac{BN-CM}{BG-CE}G.
\end{eqnarray}
This system defines discrete deformations of the structure constants.

\bigskip

\section{\textbf{Nilpotent DDA}}

For the \textbf{nilpotent} \textbf{DDA }$\mathit{L} _{3}$, in virtue of the
defining relations (32), one has

\begin{equation}
\lbrack p_{1},\varphi (x)]=0,\quad [p_{2},\varphi (x)]=\frac{\partial
\varphi }{\partial x}\cdot p_{1}
\end{equation}
or

\begin{equation}
\lbrack p_{j},\varphi (x)]=\frac{\partial \varphi }{\partial x}\cdot
\sum_{k=1}^{2}a_{jk}p_{k}
\end{equation}
where $a_{21}=1,a_{11}=a_{12}=a_{22}=0$. Using (71), one gets the following
CS

\begin{eqnarray}
\sum_{q=1}^{2}a_{lq}\sum_{m=0}^{2}C_{qm}^{n}\frac{\partial C_{jk}^{m}}{%
\partial x}-\sum_{q=1}^{2}a_{jq}\sum_{m=0}^{2}C_{qm}^{n}\frac{\partial
C_{kl}^{m}}{\partial x}+  \notag \\
+\sum_{m=0}^{2}(C_{jk}^{m}C_{lm}^{n}-C_{kl}^{m}C_{jm}^{n})=0,\quad
j,k,l,n=0,1,2.
\end{eqnarray}
In the matrix form it is

\begin{equation}
C_{1}\frac{\partial C_{1}}{\partial x}=[C_{1},C_{2}].
\end{equation}
For invertible matrix $C_{1}$

\begin{equation}
\frac{\partial C_{1}}{\partial x}=C_{1}^{-1}[C_{1},C_{2}].
\end{equation}
This system of ODEs has three independent first integrals

\begin{equation}
I_{1}=trC_{1},\quad I_{2}=\frac{1}{2}tr(C_{1})^{2},\quad I_{3}=\frac{1}{3}%
tr(C_{1})^{3}.
\end{equation}
and equivalent to the compatibility condition for the linear system

\begin{eqnarray}
C_{1}\Phi =\lambda \Phi ,  \notag \\
C_{1}\frac{\partial \Phi }{\partial x}+C_{2}\Phi =0.
\end{eqnarray}
So, as in the previous section the CS (73) describes iso-spectral
deformations of the matrix $C_{1}$. This CS governs deformations generated
by $\mathit{L} _{3}$.

For the two-dimensional algebra \textit{A} without unite element the CS is
given by equation (73) with the matrices (48). First integrals in this case
are $I_{1}=B+G,I_{2}=\frac{1}{2}(B^{2}+G^{2}+2CE)$ and $\det C_{1}=\frac{1}{2%
}I_{1}^{2}-I_{2}.$ Since $\det C_{1}$ is a constant on the solutions of the
system, then at $\det C_{1}\neq 0$ one can always introduce the variable y
defined by $x=$ $y\det C_{1}$ such that CS (74) takes the form

\begin{eqnarray}
B^{\prime} =EBG+ENC-GMC-CE^{2},  \notag \\
E^{\prime} =GBM+GEN-ECM-MG^{2},  \notag \\
C^{\prime} =BCE+BG^{2}+MC^{2}-CEG-BNC-GB^{2},  \notag \\
G^{\prime} =CMG+CE^{2}-CEN-BGE
\end{eqnarray}
where $B^{\prime} =\frac{\partial B}{\partial y}$ etc and M, N are arbitrary
functions. At $\det C_{1}=BG-CE=1$ this system becomes

\begin{eqnarray}
B^{\prime} =E+C(EN-GM),  \notag \\
E^{\prime} =M+G(EN-GM),  \notag \\
C^{\prime} =G-B+C(MC-BN),  \notag \\
G^{\prime} =-E-C(EN-GM).
\end{eqnarray}

Chosing M=N=0, one gets

\begin{equation}
B^{\prime} =E,\quad E^{\prime} =0,\quad C^{\prime} =G-B,\quad G^{\prime} =-E.
\end{equation}
The solution of this system is

\begin{equation}
E=\alpha,\quad B=\alpha y+\beta ,\quad G=-\alpha y+\gamma ,\quad
C=-y^{2}+(\gamma -\beta )y+\delta
\end{equation}
where $\alpha ,\beta ,\gamma ,\delta $ are arbitrary constants subject the
constraint $\beta \gamma -\alpha \delta =1$. First integrals for this
solution are $I_{1}=\beta +\gamma ,I_{2}=\frac{1}{2}(\beta ^{2}+\gamma
^{2}+2\alpha \delta ).$

With the choice M=0, N=1 and under the constraint $I_{1}=B+G=0$ the system
(77) takes the form

\begin{equation}
B^{\prime} =(1+C)E,\quad E^{\prime} =-BE,\quad C^{\prime} =-(2+C)B.
\end{equation}
This system can be written as a single equation in the different equivalent
forms. One of them is

\begin{equation}
(E^{\prime })^{2}+\alpha E^{4}-2E^{3}+E^{2}=0
\end{equation}
where $\alpha $ is an arbitrary constant and

\begin{equation}
B^{2}=-1-\alpha E^{2}+2E,\quad C=\alpha E-2,\quad G=-B.
\end{equation}
The second integral is equal to -1.

Solutions of equation (82) can be expressed through the elliptic integrals.
Solution of (82) and the formulae (83) define deformations of the structure
constants driven by DDA $\mathit{L} _{3}$.

\bigskip

\section{\textbf{Solvable DDAs}}

1. For the \textbf{solvable DDA }$\text{L} _{4}$ the relations (32) imply
that

\begin{equation}
\lbrack p_{j},\varphi (x)]=(T-1)\varphi (x)p_{j},\quad j=1,2
\end{equation}
where $\varphi (x)$ is an arbitrary function and $T$ is the shift operator $%
T\varphi (x)=\varphi (x+1)$. With the use of (84) one arrives at the
following CS

\begin{equation}
C_{1}TC_{2}=C_{2}TC_{1}.
\end{equation}
For nondegenerate matrix $C_{1}$ equation (85) is equivalent to the equation
$T(C_{2}C_{1}^{-1})=C_{1}^{-1}C_{2}$ or

\begin{equation}
TU=C_{1}^{-1}UC_{1}
\end{equation}
where $U\doteqdot C_{2}C_{1}^{-1}$. Using this form of the CS, one promptly
concludes that the CS (85) has three independent first integrals

\begin{equation}
I_{1}=tr(C_{2}C_{1}^{-1}),\quad I_{2}=\frac{1}{2}tr(C_{2}C_{1}^{-1})^{2},%
\quad I_{3}=\frac{1}{3}tr(C_{2}C_{1}^{-1})^{3}
\end{equation}
and is representable as the commutativity condition for the linear system

\begin{eqnarray}
\Phi C_{2}C_{1}^{-1} =\lambda \Phi ,  \notag \\
T\Phi =\Phi C_{1}.
\end{eqnarray}

For the two-dimensional algebra \textit{A} one has the CS (85) with the
matrices (48). It is the system of four equations for six functions

\begin{eqnarray}
BTE+ETG =ETB+MTC,  \notag \\
BTM+ETN =ETE+MTG,  \notag \\
CTE+GTG =GTB+NTC,  \notag \\
CTM+GTN =GTE+NTG.
\end{eqnarray}

Chosing B and C as free functions and assuming that BG-CE$\neq 0$, one can
easily resolve (89) with respect to TE,TG,TM,TN. For instance, with B=C=1
one gets the following four-dimensional mapping

\begin{eqnarray}
TE =M-E\frac{M-N}{E-G},\quad TG=1+\frac{M-N}{E-G},  \notag \\
TM =N+(N-G)\frac{M-N}{E-G}-G\left( \frac{M-N}{E-G}\right) ^{2}, \\
TN =M+(1-E)\frac{M-N}{E-G}+\left( \frac{M-N}{E-G}\right) ^{2}.  \notag
\end{eqnarray}

\ 2. In a similar manner one finds the CS associated with the \textbf{%
solvable DDA }$\text{L} _{5}$ . Since in this case

\begin{equation}
\lbrack p_{1},\varphi (x)]=(T-1)\varphi (x)p_{1},\quad [p_{2},\varphi
(x)]=(T^{-1}-1)\varphi (x)p_{2}
\end{equation}
the CS takes the form

\begin{equation}
C_{1}TC_{2}=C_{2}T^{-1}C_{1}.
\end{equation}
For nondegenerate $C_{2}$ it is equivalent to

\begin{equation}
TV=C_{2}VC_{2}^{-1}
\end{equation}
where $V\doteqdot T^{-1}C_{1}\cdot C_{2}$. Similar to the previous case the
CS has three first integrals

\begin{equation}
I_{1}=tr(C_{1}TC_{2}),\quad I_{2}=\frac{1}{2}tr(C_{1}TC_{2})^{2},\quad I_{3}=%
\frac{1}{3}tr(C_{1}TC_{2})^{3}
\end{equation}
and is equivalent to the compatibility condition for the linear system

\begin{eqnarray}
(T^{-1}C_{1})C_{2}\Phi =\lambda \Phi ,  \notag \\
T\Phi =C_{2}\Phi .
\end{eqnarray}
Note that the CS (92) is of the form (22) with $T_{1}=T,T_{2}=T^{-1}$. Thus,
the deformations generated by $\mathit{L} _{5}$ \ can be considered as the
reductions of the discrete deformations (22) under the constraint $%
T_{1}T_{2}C_{jk}^{n}=C_{jk}^{n}$.

\ A class of solutions of the CS (92) is given by

\begin{equation}
C_{j}=g^{-1}T_{j}g
\end{equation}
where g is $3\times 3$ matrix and $T_{0}=1,T_{1}=T,T_{2}=T^{-1}$. Since $%
C_{jk}^{n}=C_{kj}^{n}$ one has $T_{j}g_{l}^{m}=T_{l}g_{j}^{m}$ and hence $%
g_{j}^{m}=T_{j}\Phi ^{m}$ where $\Phi ^{0},\Phi ^{1},\Phi ^{2}$ are
arbitrary functions. So, this subclass of deformations are defined by three
arbitrary functions.

To describe the iso-associative deformations for which $%
C_{1}(x)C_{2}(x)=C_{2}(x)C_{1}(x)$ for all x these functions should obey the
systems of equations

\begin{equation}
\sum_{l,t=0}^{2}T_{j}T_{t}\Phi ^{n}\cdot (g^{-1})_{l}^{t}\cdot
T_{k}T_{m}\Phi ^{l}=\sum_{l,t=0}^{2}T_{k}T_{t}\Phi ^{n}\cdot
(g^{-1})_{l}^{t}\cdot T_{j}T_{m}\Phi ^{l},  \notag \\
\quad j,k,n,m=0,1,2.
\end{equation}
It is a version of the discrete oriented associativity equation.

\bigskip

\textbf{References}

\bigskip

1. Gerstenhaber M., On the deformation of rings and algebras, Ann. Math.,
\textbf{79}, 59-103 \ (1964).

2. Gerstenhaber M., On the deformation of rings and algebras. II, Ann.
Math., \textbf{84}, 1-19 (1966).

3. Witten E., On the structure of topological phase of two-dimensional
gravity, Nucl. Phys., \textbf{B 340}, 281-332 (1990).

4. Dijkgraaf R., Verlinde H. and Verlinde E.,\ Topological strings in $d<1$,
Nucl. Phys., \textbf{\ B 352,} 59-86 (1991).

5. Dubrovin B., Integrable systems in topological field theory, Nucl. Phys.,
\textbf{B 379}, 627-689 (1992).

6. Dubrovin B., Geometry of 2D topological field theories, Lecture Notes in
Math., \textbf{1620}, 120-348 , Springer, Berlin (1966).

7. Hertling C. and Manin Y.I., Weak Frobenius manifolds, Int. Math. Res.
Notices, \textbf{6}, 277-286 (1999).

8. Manin Y.I., Frobenius manifolds, quantum cohomology and moduli spaces,
AMS, Providence, (1999).

9. Hertling C.and Marcoli M., (Eds),\ Frobenius manifolds,quantum cohomology
\ and singularities, Aspects of Math., \textbf{E36}, Friedr. Vieweg \& Sohn,
\ Wiesbaden (2004).

10. Konopelchenko B.G. and \ Magri F., Coisotropic deformations of
associative algebras and dispersionless integrable hierarchies, Commun.
Math. Phys., \textbf{274}, 627-658 (2007)

11. Konopelchenko B.G. and Magri F., Dispersionless integrable equations as
coisotropic deformations: extensions and reductions, Theor. Math. Phys.,
\textbf{151}, 803-819 (2007).

12. Konopelchenko B.G., Quantum deformations of associative algebras and
integrable systems, J. Phys. A: Math. Theor., \textbf{42}, 095201 (2009);
arXiv:0802.3022 (2008).

13. Konopelchenko B.G., \ Discrete integrable systems and deformations of
associative algebras, arXiv:0904.2284 (2009).

14. Konopechenko B.G. and Magri F., Yano manifolds, F-manifolds and
integrable systems, to appear.

15. Van del Waerden \ B.L., Algebra, Springer-Verlag, 1971.

16. Diedonne J., Les determinants sur un corps noncommutatiff, Bul. Soc.
Math. France, \textbf{71}, 27-45 (1943).

17. Gelfand I.M. and Retakh V.S., Determinants of matrices over
noncommutative rings, Funkt. Anal. Appl., \textbf{25}, 13-25
(1991).

18. Bourbaki N., Groupes et algebres de Lie, Hermann, Paris, (1972).

19. Dirac P.A.M., Lectures on quantum mechanics, Yeshiva Univ., New York
(1964).

20. Losev A. and Manin Y.I., \ Extended modular operads, in: Frobenius
manifolds, quantum cohomology and singularities, ( Eds. C. Hertling and M.
Marcoli), Aspects of Math., \textbf{E36}, 181-211(2004).

21. Yano K. and Ako M., On certain operators associated with tensor fileds,
Kodai Math.Sem. Rep., \ \textbf{20} , 414-436 (1968).

22. Novikov S.P., Manakov S.V., Pitaevski L.P. and Zakharov V.E., Theory of
solitons. The inverse problem method, Nauka, (1980); Plenum, New York (1984).

23. Ablowitz M.J. and Segur H., Solitons and inverse scattering transform ,
SIAM, Philadelphia (1981).

24. Baker H.F., Abelian functions, Cambridge Univ. Press (1995).

25. Chazy J., Sur les equations differentielles du troisieme ordre et
d'ordre superieur dont l'integrale generale a ses points critiques fixes,
Acta Math., \textbf{34}, 317-385 (1911).

26. Bureau F., Differential equations with fixed critical points, Annali
Matemat., \textbf{66}, 1-116 (1964).

27. Cosgrove C.M., Chazy classes IX-XI of third order differential
equations, Stud.Applied Math., \textbf{104}, 171-228 (2000).

\end{document}